\documentclass[a4paper,12pt]{amsart}

\usepackage{amssymb}
\usepackage{latexsym} 
\usepackage{amsfonts} 
\usepackage{amsmath}
\usepackage{eucal} 
\usepackage{bm} 
\usepackage{bbm} 
\usepackage{graphicx} 
\usepackage[english]{varioref} 
\usepackage[nice]{nicefrac} 
\usepackage[all]{xy}
\usepackage{amsthm}

\def\i{^{-1}}
\def\ge{\geqslant}
\def\le{\leqslant}
\def\<{\langle} 
\def\>{\rangle}

\def\a{\alpha}

\def\e{\epsilon}

\def\s{\sigma}

\def\k{\kappa}
\def\l{\lambda}

\def\Om{\Omega}

\def\ZZ{\mathbb Z}
\def\NN{\mathbb N}
\def\QQ{\mathbb Q}
\def\FF{\mathbb F}
\def\RR{\mathbb R}

\def\kk{\mathbf k}

\def\co{\mathcal O}

\def\tH{\tilde H}

\def\tW{\tilde W}

\def\subset{\subseteq}

\theoremstyle{plain}
\newtheorem{thm}{Theorem}[section] 
\newtheorem*{thm*}{Theorem}

 \newtheorem{cor}[thm]{Corollary}
  \newtheorem{conj}[thm]{Conjecture}

\theoremstyle{definition}

\theoremstyle{remark}

\newtheorem*{claim*}{Claim}

\begin{document}

\author{Xuhua He}
\address{Department of Mathematics and Institute for Advanced Study, The Hong Kong University of Science and Technology, Clear Water Bay, Kowloon, Hong Kong}
\email{maxhhe@ust.hk}
\thanks{The author is partially supported by HKRGC grant 602011.}
\title[Note on affine Deligne-Lusztig varieties]{Note on affine Deligne-Lusztig varieties}

\begin{abstract}
This note is based on my talk at ICCM 2013, Taipei. We give an exposition of the group-theoretic method and recent results on the questions of non-emptiness and dimension of affine Deligne-Lusztig varieties in affine flag varieties. 
\end{abstract}

\maketitle

\section*{Introduction}
Affine Deligne-Lusztig varieties was first introduced by Rapoport in \cite{Ra}, generalizing Deligne and Lusztig's classical
construction. Understanding the emptiness/non-emptiness pattern and dimension of affine Deligne-Lusztig varieties is fundamental in the study of reduction of Shimura varieties with parahoric level structures.

In this note, we will discuss the new group-theoretic method to study affine Deligne-Lusztig varieties in affine flag varieties and answer the above questions in terms of class polynomials of affine Hecke algebras.  The two key ingredients are the Deligne-Lusztig reduction \cite{DL} and the combinatorial properties of affine Weyl groups \cite{HN}. 

\section{Affine Deligne-Lusztig varieties}

\subsection{} Let $\FF_q$ be a finite field with $q$ elements and $\kk$ be an algebraic closure of $\FF_q$. We consider the field $L=\kk((\e))$ and its subfield $F=\FF_q((\e))$. The Frobenius automorphism $\s$ of $\kk/\FF_q$ can be extended in the usual way to an automorphism of $L/F$ such that $\s(\sum a_n \e^n)=\sum \s(a_n) \e^n$. 

Let $G$ be a split connected reductive group over $\FF_q$ and let $I$ be a $\s$-stable Iwahori subgroup of $G(L)$.  The affine Deligne-Lusztig variety associated with $w$ in the extended affine Weyl group $\tW\cong I \backslash G(L)/I$ and $b \in G(L)$ is defined to be $$X_{w}(b)=\{g I \in G(L)/I; g \i b \s(g) \in I w I\}.$$ Here we embed $\tW$ set-theoretically into $G(L)$. 

This is a locally closed sub-ind scheme of the affine flag variety $G(L)/I$. It is a finite dimensional $\kk$-scheme, locally of finite type over $\kk$.\footnote{One may consider any connected reductive group over $F$ that splits over a tamely ramified extension of $L$. As discussed for example in \cite[Section 6 \& 7]{He99}, the questions we'll discuss in this note can be reduced to quasi-split unramified groups. For simplicity, we only consider split groups in this note.}

\subsection{} To illustrate some difficulties in the study of affine Deligne-Lusztig varieties, let us begin with (classical) Deligne-Lusztig varieties. Let $B$ be a Borel subgroup of $G$ defined over $\FF_q$, and $W$ be the associate finite Weyl group. For any $w \in W$ and $b \in G(\kk)$, the classical Deligne-Lusztig variety $X_{w, b}$ is defined to be $$X_{w, b}=\{g B \in G(\kk)/B(\kk); g \i b \s(g) \in B(\kk) w B(\kk)\}.$$ 

Note that in both the classical case and the affine cases, if $b$ and $b'$ are in the same $\s$-conjugacy class, then the associated classical/affine Deligne-Lusztig varieties are isomorphic. In other words, the classical/affine Deligne-Lusztig varieties depends on the finite/affine Weyl group element $w$ and the $\s$-conjugacy class $[b]$. 

However, in the classical case, $G(\kk)$ is a single $\s$-conjugacy class. Therefore, we may omit $b$ in the definition and simply write $X_w=\{g B \in G(\kk)/B(\kk); g \i \s(g) \in B(\kk) w B(\kk)\}$. This is what appears in the literature. The classical Deligne-Lusztig variety $X_w$ is always nonempty and is a smooth variety of dimension $\ell(w)$. 

The loop group $G(L)$, on the other hand, contains infinitely many $\s$-conjugacy classes. The extra parameter $[b]$ makes the study of affine Deligne-Lusztig varieties much harder than that of the classical ones. 

From a different point of view, the two parameters $w$ and $[b]$ in the definition of affine Deligne-Lusztig varieties appear naturally. There are two important stratifications on the special fiber of a Shimura variety with Iwahori level structure: one is the Kottwitz-Rapoport stratification whose strata are indexed by specific elements $w \in \tW$; the other is the Newton stratification whose strata are indexed by specific $\s$-conjugacy classes $[b] \subset G(L)$. There is a close relation between the affine Deligne-Lusztig variety $X_{w}(b)$ and the intersection of the Newton stratum associated with $[b]$ with the Kottwitz-Rapoport stratum associated with $w$ (see \cite{Ha} and \cite{VW}). 

\section{Conjugacy classes of $\tW$}

\subsection{}\label{2.1} To understand the affine Deligne-Lusztig varieties $X_w(b)$, one needs to understand the $\s$-conjugacy classes of $G(L)$ and its relation with the extended affine Weyl group $\tW$. We will relate the conjugacy classes of $\tW$ with the $\s$-conjugacy classes of $G(L)$. Based on the decomposition $G(L)=\sqcup_{w \in \tW} I w I$, this sounds plausible. However, the naive map sending $I w I$ to the $\s$-conjugacy class of $w$ does not work. The reasons are as follows:

(1) Given two elements $w$ and $w'$ in the same conjugacy class of $\tW$, the set $I w I$ and $I w' I$ may not be $\s$-conjugated to each other. 

(2) One $\s$-conjugacy class of $G(L)$ may contain elements of different conjugacy classes in $\tW$. For example, in the classical case, $G(\kk)$ is a single $\s$-conjugacy class and $W$ contains several conjugacy classes. 

The ideas to overcome the difficulties is to use the minimal length elements in the conjugacy classes of $\tW$. We will recall some properties of the minimal length elements in this section and discuss the applications to $\s$-conjugacy classes of $G(L)$ in the next section. 

\subsection{} We have the semidirect product $$\tW=P \rtimes W=\{t^\l w; \l \in P, w \in W\},$$ where $P$ is the coweight lattice of a split maximal torus $T$ of $G$. Let $Q$ be the coroot lattice of $T$. Then $W_a=Q \rtimes W \subset \tW$ is an affine Weyl group. The length function and Bruhat order on $W_a$ extends in a natural way to $\tW$ and $\tW=W_a \rtimes \Om$, where $$\Om=\{w \in \tW; \ell(w)=0\} \cong P/Q.$$ 

Let $\k: \tW \to \tW/W_a$ be the natural projection. We call $\k$ the Kottwitz map. 

Let $P_{\QQ}=P \otimes_\ZZ \QQ$ and $P_\QQ/W$ be the quotient of $P_\QQ$ by the natural action of $W$. We may identify $P_\QQ/W$ with $P_{\QQ, +}$, the set of dominant rational coweights. For any $w \in \tW$, $w^{n_0}=t^\l$ for some $\l \in P$, where $n_0=\sharp(W)$. Let $\nu_{w}$ be the unique element in $P_{\QQ, +}$ that lies in the $W$-orbit of $\l/n_0$. We call the map $$\tW \to P_{\QQ, +}, \qquad w \mapsto \nu_w$$ the Newton map. 

Define $f: \tW \to \tW/W_a \times P_{\QQ, +}$ by $w \mapsto (\k(w), \nu_w)$. It is constant on each conjugacy class of $\tW$. We denote the image of $f$ by $B(\tW)$. Note that in general, a fiber of $f$ contains more than one conjugacy class of $\tW$. 

\subsection{}\label{straight} We call an element $w \in \tW$ straight if $\ell(w)=\<\nu_w, 2 \rho\>$, where $\rho$ is the half sum of all the positive roots. It is easy to see that $w$ is straight if and only if $\ell(w^n)=n \ell(w)$ for all $n \in \NN$. We call a conjugacy class of $\tW$ straight if it contains some straight element. Note that a straight element is automatically of minimal length in its conjugacy class. 

By \cite[Proposition 3.2]{HN}, a conjugacy class is straight if and only if it contains a length-zero element in the extended affine Weyl group $\tilde W_M$ for some standard Levi subgroup $M$ of $G$. 

By \cite[Theorem 3.3]{HN}, the map $f: \tW \to \tW/W_a \times P_{\QQ, +}$ induces a bijection from the set of straight conjugacy classes to $B(\tW)$.\footnote{In \cite[Section 3]{HN}, we assume that $G$ is adjoint. However, \cite[Theorem 3.3]{HN} holds for any reductive group and the proof is the same as in loc. cit.}

\subsection{} Now we discuss the minimal length elements in a conjugacy class $\co$ of $\tW$. The following result is obtained in \cite{HN}, motivated by previous results of Geck and Pfeiffer \cite{GP93} for finite Coxeter groups. 

\begin{thm}\label{min}
Let $\co$ be a conjugacy class of $\tW$. Then 

(1) Each element of $\co$ can be brought to a minimal length element by conjugation by simple reflections which reduce or keep constant the length. 

(2) Any two minimal length elements in $\co$ are conjugate in the associated Braid group. 

(3) If moreover, $\co$ is straight, then any two minimal length elements in $\co$ are conjugate by cyclic shifts. 
\end{thm}

Let $\co$ be a straight conjugacy class and $\co'$ be another conjugacy class such that $f(\co)=f(\co')$. Then the minimal length elements in $\co'$ are related to the straight elements in $\co$ in the sense of \cite[Theorem 3.4]{HN}. This property is used to overcome the difficulty $\S$\ref{2.1} (2). 

\section{$\s$-conjugacy classes of $G(L)$}

\subsection{}\label{3.1} Recall that $G(L)=\sqcup_{w \in \tW} I w I$. Let $w \in \tW$ and $s \in \tW$ be a simple reflection. Then \[I s I w I=\begin{cases} I s w I, & \text{ if } s w>w; \\ I s w I \sqcup I w I, & \text{ if } s w<w.\end{cases}\] \[I w I s I=\begin{cases} I w s I, & \text{ if } w s>w; \\ I w s I \sqcup I w I, & \text{ if } w s<w.\end{cases}\]

Therefore, \[G(L) \cdot_\s I w I=\begin{cases} G(L) \cdot_\s I s w s I, & \text{ if } \ell(s w s)=\ell(w); \\ G(L) \cdot_\s I s w s I \cup G(L) \cdot_\s I s w I, & \text{ if } \ell(s w s)<\ell(w). \end{cases}\] Here $\cdot_\s$ denotes the $\s$-conjugation action. 

This equality, together with Theorem \ref{min}, gives a reduction method in the study of $G(L) \cdot_\s I w I$ and allows us to reduce the general case to the case where $w$ of minimal length in its conjugacy class. This is used to overcome the difficulty $\S$\ref{2.1} (1). 

We have the following parameterisation of $\s$-conjugacy classes of $G(L)$. 

\begin{thm}\label{para}
There is a canonical bijection between 

(a) The set of $\s$-conjugacy classes of $G(L)$; 

(b) The set of straight conjugacy classes of $\tW$;

(c) The image of $f: \tW \to \tW/W_a \times P_{\QQ, +}$.
\end{thm}

Here the bijection between (a) and (c) follows from Kottwitz's classification of $\s$-conjugacy classes \cite{Ko97} together with the fact that any $\s$-conjugacy class is represented by an element in $\tW$ \cite[Corollary 7.2.2]{GHKR}. The bijection between (b) and (c) is discussed in $\S$\ref{straight}. The bijection between (a) and (b) (and a new proof of the classification of $\s$-conjugacy classes) is obtained in \cite[Section 3]{He99} using the strategy in $\S$\ref{3.1}. 

\subsection{} Now we introduce partial orders on the three sets in Theorem \ref{para}. 

Let $C, C'$ be $\s$-conjugacy classes of $G(L)$. We write $C \le C'$ if $C$ is contained in the closure of $C'$. 

Let $\co, \co'$ be straight conjugacy classes of $\tW$. We write $\co \le \co'$ if for some $w'$ of minimal length in $\co'$, there exists $w$ of minimal length in $\co$ such that $w \le w'$ with respect to the Bruhat order in $\tW$. It is proved in \cite[$\S$4.7]{HeMin} and \cite[Corollary 7.5]{He5} that $\le$ is a partial order on the set of straight conjugacy classes. 

Let $(k, \nu), (k', \nu') \in \tW/W_a \times P_{\QQ, +}$. We write $(k, \nu) \le (k', \nu')$ if $k=k'$ and $\nu'-\nu \in \sum_{\a} \RR_{\ge 0} \a$, where $\a$ runs over all the simple roots. This partial order is studied in detail by Chai in \cite{C}. 

\begin{thm}
Let $C, C'$ be $\s$-conjugacy classes of $G(L)$ and $\co, \co'$ the corresponding straight conjugacy classes in $\tW$. Then the following conditions are equivalent:

(1) $C \le C'$;

(2) $\co \le \co'$;

(3) $f(\co) \le f(\co')$. 
\end{thm}

Here $(1) \Rightarrow (3)$ is proved by Rapoport and Richartz in \cite[Theorem 3.6]{RR}, $(3) \Rightarrow (1)$ is proved by Viehmann in \cite[Theorem 2]{Vi} and the equivalence between (1) and (2) is proved in \cite[Section 11]{He5} and \cite{He00}. It is easy to show that $(2) \Rightarrow (3)$. It is interesting to give a direct (combinatorial proof) that $(3) \Rightarrow (2)$.\footnote{In fact, $(3) \Leftrightarrow (1) \Rightarrow (2)$ holds for any tamely ramified group. We expect that $(2) \Rightarrow (1)$ also holds for any tamely ramified group. As explained in the beginning of this note, it suffices to prove it for quasi-split unramified groups.}

\section{``Dimension=degree'' Theorem}

\subsection{} We first recall the reduction method of Deligne and Lusztig \cite[Theorem 1.6]{DL}. 

\begin{thm}
Let $w\in \tW$, and let $s$ be a simple affine reflection.
\begin{enumerate}
\item 
If $\ell(sws) = \ell(w)$, then there exists a universal homeomorphism $X_w(b) \rightarrow X_{sws}(b)$.
\item
If $\ell(sws) = \ell(w)-2$, then $X_w(b)$ can be written as a disjoint union $X_w(b) = X_1 \sqcup X_2$ where $X_1$ is closed and $X_2$ is open, and such that there exist morphisms $X_1\to X_{sws}(b)$ and $X_2\to X_{sw}(b)$ which are compositions of a Zariski-locally trivial fiber bundle with one-dimensional fibers and a universal homeomorphism.
\end{enumerate}
\end{thm}

The reduction method, together with Theorem \ref{min}, in principle, reduce the study of $X_w(b)$ to the case where $w$ is of minimal length in its conjugacy class. The latter one, is studied in detail in \cite{HL} and \cite{He99}. In particular, 

\begin{thm}
Let $w \in \tW$ be an element of minimal length in its conjugacy class and $b \in G(L)$. 

(1) If $b$ and $w$ are not in the same $\s$-conjugacy class of $G(L)$, then $X_w(b)=\emptyset$. 

(2) If $b$ and $w$ are in the same $\s$-conjugacy class of $G(L)$, then $\dim X_w(b)=\ell(w)-\<\nu_w, 2 \rho\>$.
\end{thm}

The emptiness/nonemptiness pattern and dimension formula for affine Deligne-Lusztig varieties is obtained if we can keep track of the reduction step from an arbitrary element to a minimal length element. This is accomplished via the class polynomials of affine Hecke algebras. 

\subsection{} Let $\tH$ be the Hecke algebra associated with $\tW$, i.e., $\tH$ is the associated $A=\ZZ[v, v \i]$-algebra with basis $T_{w}$ for $w \in \tW$ and multiplication is given by 
\begin{gather*} T_{w} T_{w'}=T_{w w'}, \quad \text{ if } \ell(w)+\ell(w')=\ell(w w'); \\ (T_s-v)(T_s+v \i)=0, \quad \text{ for any simple reflection } s. \end{gather*}

It is proved in \cite{HN} that for any $w \in \tW$
and conjugacy class $\co$ of $\tW$, there exists a unique polynomial $f_{w, \co} \in \NN[v-v \i]$
such that for any finite dimensional representation $V$ of $\tH$, 
\[Tr(T_w, V)=\sum_{\co} f_{w, \co} Tr(T_{w_\co}, V), \]
where $w_\co$ is a minimal length element in $\co$.

\subsection{} Now  we state the main result in \cite{He99}, which relates the dimension of affine Delgine-Lusztig varieties with the degree of the class polynomials. 

\begin{thm}\label{class3}
Let $b \in G(L)$ and $w \in \tW$. Then 

\[\dim (X_{w}(b))=\max_{\co} \frac{1}{2}(\ell(w)+\ell(\co)+\deg(f_{w, \co}))-\<\bar \nu_b, 2\rho\>,\] here $\co$ runs over conjugacy classes of $\tW$ with $f(\co)=f(b)$ and $\ell(\co)$ is the length of any minimal length element in $\co$.
\end{thm}

Here we use the convention that the dimension of an empty variety and the degree of a zero polynomial are both $-\infty$. So in particular, $X_w(b) \neq \emptyset$ if and only if $f_{w, \co} \neq 0$ for some conjugacy class $\co$ of $\tW$ with $f(\co)=f(b)$. 

\section{Affine Deligne-Lusztig varieties for basic $b$}

\subsection{} We say that an element $b \in G(L)$ is basic if $\<\nu_b, \a\>=0$ for any root $\a$ of $G$. For basic $b$, we are able to give a more explicit description of the emptiness/nonemptiness pattern and dimension formula, as conjectured by  G\"ortz, Haines, Kottwitz and Reuman  in \cite[Conjecture 1.1.1 \& Conjecture 1.1.3]{GHKR}. We first discuss the emptiness/nonemptiness pattern. It is given in terms of the $P$-alcoves introduced in \cite[Definition 2.1.1]{GHKR}. 

Let $P=M N$ be a semistandard parabolic subgroup of $G$ and $w \in \tW$. We say that $w$ is a $P$-alcove element if $w \in \tW_M$ and $N(L) \cap w I w \i \subset N(L) \cap I$.  

\begin{thm}\label{nonempty}
Let $w \in \tW$ and $b \in G(L)$ be a basic element. Then $X_w(b) \neq \emptyset$ if and only if for any semistandard parabolic subgroup $P=M N$ such that $w$ is a $P$-alcove element, $\k_M(b) \in \k_M([b] \cap M(L))$. 
\end{thm}

The ``only if'' part is proved in \cite[Theorem 1.1.2]{GHKR} as a consequence of the Hodge-Newton decomposition for affine Deligne-Lusztig varieties \cite[Theorem 1.1.4]{GHKR}. The ``if'' part is proved in \cite[Theorem A]{GHN} by showing that the notion of $P$-alcoves is compatible with the Deligne-Lusztig reduction. An algebraic proof of the analogy of Hodge-Newton decomposition for affine Hecke algebras and a new proof of the ``only if'' part of the Theorem \ref{nonempty} is obtained in \cite{HN2}.

\subsection{} 

Define

\begin{itemize}
\item $\eta_1: \tilde W=P \rtimes W  \to W$, the projection map.

\item $\eta_2: \tilde W \to W$ such that $\eta_2(w) \i w$ lies in the dominant Weyl chamber.

\item $\eta(w)=\eta_2(w) \i \eta_1(w) \eta_2(w)$.
\end{itemize}

As a consequence of Theorem \ref{nonempty}, there is a simpler criterion for emptiness/nonemptiness if $w$ lies in the shrunken Weyl chamber. 

\begin{cor}
Assume that the Dynkin diagram of $G$ is connected. Let $b \in G(L)$ be a basic element and $w \in \tW$ lies in the shrunken Weyl chamber. Then $X_w(b) \neq \emptyset$ if and only if $\k(b)=\k(w)$ and $\eta(w)$ is not in any proper parabolic subgroup of $W$.  
\end{cor}

\subsection{} The following dimension formula is conjectured in \cite[Conjecture 1.1.3]{GHKR} and proved in \cite[Theorem 12.1]{He99}. 

\begin{thm}
Let $b \in G(L)$ be a basic element and $w \in \tW$ lie in the shrunken Weyl chamber. If $X_w(b) \neq \emptyset$, then $$\dim X_w(b)=\frac{1}{2}(\ell(w)+\ell(\eta(w))-\text{def}(b)),$$ where $\text{def}(b)$ is the defect of $b$. 
\end{thm}

Here the lower bound is obtained by constructing an explicit sequence from $w$ to a minimal length element $w'$ in some conjugacy class of $\tW$ (in most cases, $w$ and $w'$ are not in the same conjugacy class). The upper bound is obtained by combining the ``partial conjugation'' method in \cite{HeMin} and the dimension formula for affine Deligne-Lusztig varieties in affine Grassmannian in \cite[Theorem 2.15.1]{GHKR06} and \cite[Theorem 1.1]{Vi0}. If $w$ lies in the shrunken Weyl chamber, then the lower bound and the upper bound coincide and the theorem is proved.\footnote{I learned from E. Viehmann \cite{Vi2} that her student Paul Hamacher recently proved the conjectural dimension formula for affine Deligne-Lusztig varieties in affine Grassmannian for unramified groups. Combining this result, with the proof in \cite{He99} as we outlined above, and \cite[Proposition 2.5.1]{GHN}, we are able to generalize the above theorem to any tamely ramified groups.}

It is a challenging problem to give a close formula of $\dim X_w(b)$ for $w$ in the critical strips. 

\section{Affine Deligne-Lusztig varieties for nonbasic $b$}

\subsection{} An element $w \in \tW$ can be written in a unique way as $x t^\mu y$, where $x, y \in W$ and $\mu \in P$ such that $t^\mu y$ sends the fundamental alcove to an alcove in the dominant Weyl chamber. In this case, $\eta(x t^\mu y)=y x$. 

Let $w_0$ be the longest element in $W$. For any $\mu \in P_+$, $w_0 t^\mu$ is the unique maximal element in its $W \times W$-coset. For such element, there is a complete answer for the emptiness/nonemptiness pattern and dimension formula. 

\begin{thm}
Let $b \in G(L)$ and $\mu \in P_+$. Then $X_{w_0 t^\mu}(b) \neq \emptyset$ if and only if $f(b) \le f(t^\mu)$. In this case, $$\dim X_{w_0 t^\mu}(b)=\<\mu-\nu_b, \rho\>+\ell(w_0)-\frac{1}{2}\text{def}(b).$$
\end{thm}

The affine Deligne-Lusztig variety $X_{w_0 t^\mu}(b)$ in the affine flag and the affine Deligne-Lusztig variety $X_\mu(b)$ in the affine Grassmannian are related as in \cite[Theorem 10.1]{He99}. For affine Deligne-Lusztig varieties in the affine Grassmannian, the emptiness/nonemptiness is obtained by Gashi in \cite[Theorem 1.1]{Ga} and the dimension formula is obtained by G\"ortz, Haines, Kottwitz and Reuman \cite[Theorem 2.15.1]{GHKR06} and Viehmann \cite[Theorem 1.1]{Vi0}. 

\subsection{} For other $w$, not much is known. G\"ortz, Haines, Kottwitz and Reuman made the following conjecture \cite[Conjecture 9.5.1 (b)]{GHKR} which relates the affine Deligne-Lusztig variety $X_w(b)$ with $X_w(b')$ for basic $b'$. 

\begin{conj}
Let $b \in G(L)$ and $b'$ be a basic element in $G(L)$ such that $\k(b)=\k(b')$. Then for $w \in \tW$ with sufficiently large length, $X_w(b) \neq \emptyset$ if and only if $X_w(b') \neq \emptyset$. In this case, $$\dim X_w(b)=\dim X_w(b')-\<\nu_b, \rho\>+\frac{1}{2}(\text{def}(b')-\text{def}(b)).$$
\end{conj}

Many numerical evidence for group of small rank is obtained by computer in support of this conjecture. Another evidence for split $b$ is obtained in \cite{He00}. 

\begin{thm}
Assume that the Dynkin diagram of $G$ is connected. Let $\mu \in P_+$ and $\l \in Q$ be a dominant and regular coweight. Then for any $x, y \in W$, $X_{x t^{\mu+\l} y}(t^\mu) \neq \emptyset$ if and only if $y x$ is not in any proper parabolic subgroup of $W$. In this case, $$\dim X_{x t^{\mu+\l} y}(t^\mu)=\<\l, \rho\>+\frac{1}{2}(\ell(x)+\ell(y)+\ell(y x)).$$
\end{thm}

\section*{Acknowledgment} We are grateful to M. Rapoport and E. Viehmann for their comments on this note.

\end{document}